\documentclass{elsart}

\usepackage{amsmath}
\usepackage{amssymb}
\usepackage[dvips]{graphicx}

\def\<{\scriptstyle \leftarrow}
\def\>{\scriptstyle \rightarrow}
\def\={\scriptstyle \rightleftarrows}
\def\eps{\varepsilon}

\begin{document}

\begin{frontmatter}

\title{Accuracy and convergence of \\ the backward Monte-Carlo method}

\author{Johan Carlsson\thanksref{email}}

\address{Oak Ridge National Laboratory,
P.O.~Box 2009, Oak Ridge, TN 37831--8071, USA}

\thanks[email]{E-mail: \textsf{carlssonja@ornl.gov}}

\begin{abstract}
The recently introduced backward Monte-Carlo method~[Johan Carlsson,\linebreak
\texttt{arXiv:math.NA/0010118}] is validated,
benchmarked, and compared to the conventional, forward Monte-Carlo method by
analyzing the error in the Monte-Carlo solutions to a simple model equation.
In particular, it is shown how the backward method reduces the statistical
error in the common case where the solution is of interest in only a small
part of phase space. The forward method requires binning of particles, and
linear interpolation between the bins introduces an additional error. Reducing
this error by decreasing the bin size increases the statistical error. The
backward method is not afflicted by this conflict. Finally, it is shown how
the poor time convergence can be improved for the backward method by a minor
modification of the Monte-Carlo equation of motion that governs the
stochastic particle trajectories. This scheme does not work for the
conventional, forward method. \\[4mm]

\noindent\emph{PACS:} 02.70.Lq; 02.60.Lj \\[-1mm]

\end{abstract}

\begin{keyword}
Linear parabolic partial differential equation; Monte-Carlo method;
Time convergence
\end{keyword}

\end{frontmatter}

\section{Introduction}

The simplest parabolic equation is the canonical diffusion equation,
\begin{equation}
\frac{\partial f}{\partial t} =
\frac{\partial}{\partial x} \, D \, \frac{\partial f}{\partial x} \ ,
\quad x \in \mathbb{R} \ , \quad 0 \leq t \leq T \ ,
\label{eq:canonical_diffusion}
\end{equation}
with some initial condition $f(x,0) = \Phi(x)$. The conventional, forward,
weighted Monte-Carlo solution to Eq.~(\ref{eq:canonical_diffusion}) is given
by
\begin{equation}
f_{\!\>}(x,T) = N^{-1} \sum_{i=1}^{N} \Phi{}(X_{i}^{\>}(0)) \times
\delta (x - X_{i}^{\>}(T)) \ , \label{eq:forward_solution}
\end{equation}
where the stochastic variables $X_{i}^{\>}(T)$ are found by following the
stochastic trajectories given by the forward Monte-Carlo difference equation
of motion:
\begin{equation}
X_{i}^{\>}(t + \Delta{}t) =
X_{i}^{\>}(t) + \mu \Delta{}t + \zeta \sigma \sqrt{\Delta{}t} \ ,
\quad i = 1 , \, \ldots , N \ , \label{eq:forward_motion}
\end{equation}
where $\mu = \partial D / \partial x$, $\sigma = \sqrt{2 D}$, and $\zeta$
is a zero-mean, unit-variance Gaussian random number, $\zeta \in N(0,1)$.
The points where the particles are launched, $X_{i}^{\>}(0)$, must be chosen
so that $f_{\!\>}(x,T\!=\!0)$, as given by Eq.~(\ref{eq:forward_solution}),
approximates $\Phi{}(x)$. The algorithm is illustrated in
Fig.~\ref{fig:forward} below.
\begin{figure}[!b]
\begin{center}
\includegraphics[width=120mm]{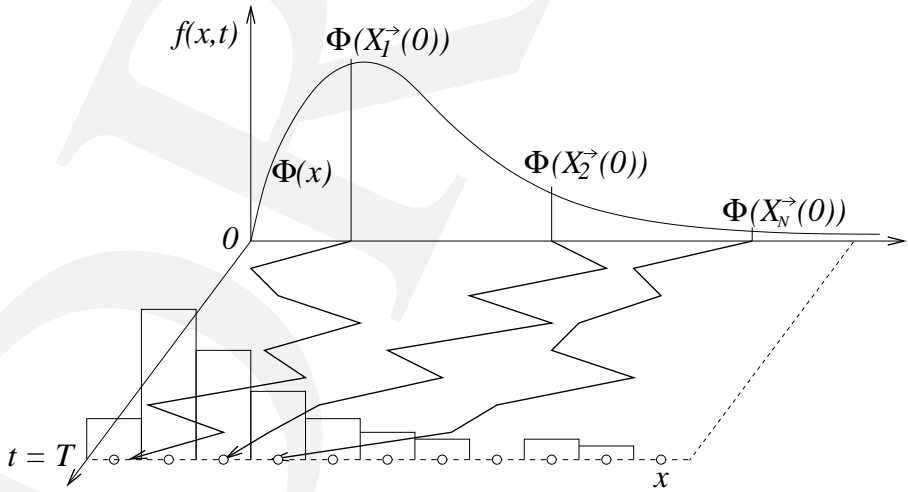}
\end{center}
\caption{The conventional, forward, weighted Monte-Carlo method.
The trajectories as given by the Monte-Carlo equation of
motion~(\ref{eq:forward_motion}).}
\label{fig:forward}
\end{figure}

The recently introduced \emph{backward} Monte-Carlo method~\cite{bmc} is
both strikingly similar and fundamentally different from its older, forward
sibling. The backward solution to Eq.~(\ref{eq:canonical_diffusion}) is
\begin{equation}
f_{\!\<}(x,T) = N^{-1} \sum_{i=1}^{N} \Phi{}(X_{i}^{\<}(0)) \ ,
\label{eq:backward_solution}
\end{equation}
where the stochastic variables $X_{i}^{\<}(0)$ are found by following the
stochastic trajectories given by the \emph{backward} Monte-Carlo difference
equation of motion:
\begin{equation}
X_{i}^{\<}(t - \Delta{}t) = X_{i}^{\<}(t) + \mu \Delta{}t +
\zeta \sigma \sqrt{\Delta{}t} \ ,
\quad X_{i}^{\<}(T) = x \, , \ i = 1 , \, \ldots , N \ .
\label{eq:backward_motion}
\end{equation}
\begin{figure}
\begin{center}
\includegraphics[width=120mm]{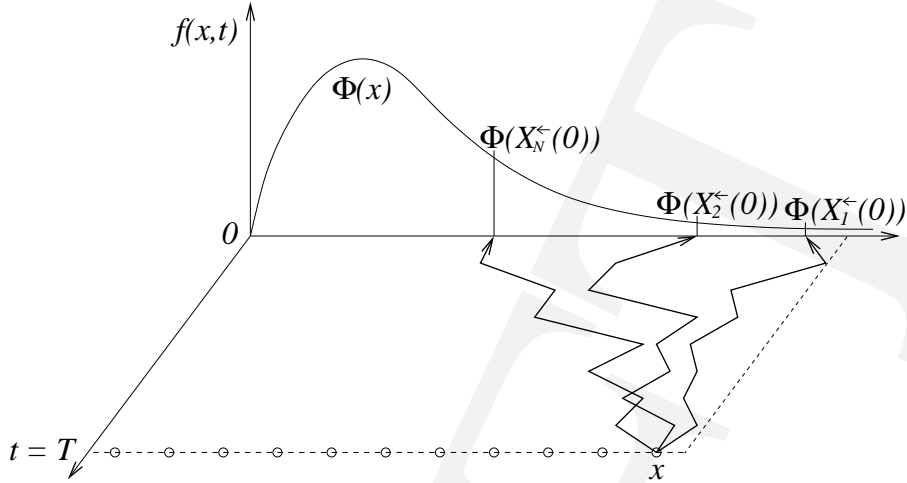}
\end{center}
\caption{The backward Monte-Carlo method. The trajectories as given by the
backward Monte-Carlo equation of motion~(\ref{eq:backward_motion}).}
\label{fig:backward}
\end{figure}

Comparing the forward and backward solutions, Eqs.~(\ref{eq:forward_solution})
and (\ref{eq:backward_solution}), respectively, the appearance of $\delta$-%
functions in the forward solution is the most striking difference. To
understand why the backward method produces a smooth solution, i.e.~without
$\delta$-functions, we have to go back to its derivation~\cite{bmc}, which
relies heavily on the Feynman-Kac formula~\cite{ito}. In fact, Eq.~%
(\ref{eq:backward_solution}), supplemented by Eq.~(\ref{eq:backward_motion}),
is just the numerical approximation of the stochastic Feynman-Kac
representation of Eq.~(\ref{eq:canonical_diffusion}). The Feynman-Kac formula
relates a parabolic equation to its \emph{naturally associated stochastic
differential equation} (SDE). This SDE governs microscopic motion going
\emph{backward} in time, and Eq.~(\ref{eq:backward_motion}) is its difference
approximation. As a consequence of time running in different directions in the
macroscopic and microscopic world, we get the natural initial condition,
$X_{i}^{\<}(T) = x$ in Eq.~(\ref{eq:backward_motion}); i.e.~we launch the
particles from the point in extended phase space where we want to know the
solution to the parabolic equation and let them work their way backward in
time to sample from the initial condition $\Phi{}(X_{i}^{\<}(0))$
[see Fig.~\ref{fig:backward} above]. As is shown in appendix~\ref{app:fmc},
the forward Monte-Carlo method can be derived by forcing the naturally
associated SDE to govern microscopic motion going \emph{forward} in time.
This coercion is responsible for both the $\delta$-functions in the forward
Monte-Carlo solution, Eq.~(\ref{eq:forward_solution}), and the
``fuzzy initial condition'',
$X_{i}^{\>}(0) \in \{X_{i}^{\>}(0) \ | \ f_{\!\>}(x,T\!=\!0) \approx
\Phi{}(x)\}$, imposed on the forward Monte-Carlo equation of motion,
Eq.~(\ref{eq:forward_motion}) [see appendix~\ref{app:fmc} for details]. \\

However, the backward algorithm was not developed for esthetic reasons.
The motivation originally came from the frustration over the extremely
inefficient use of test particles of the forward method in cases where
only the solution in a small part of phase space contributes significantly
to the physics of interest. An example of this type of problem is the heating
of fusion plasmas by cyclotron-resonance absorption on ions, where only the
high-energy tail of the solution is of relevance; but most test particles are
wasted on the almost Maxwellian bulk, which has been simulated
with e.g.~the FIDO forward Monte-Carlo code~\cite{fido}. Earlier claims~%
\cite{bmc} that the backward method is vastly superior in these cases will be
substantiated in section~\ref{sec:model}. \\

Before that, the sources of numerical error, for both the forward and the
backward method, will be identified and briefly discussed in the next section
below. In section~\ref{sec:model} we will compare the forward and backward
Monte-Carlo solutions to a simple model equation. In section~\ref{sec:higher},
it will be shown how a minor modification of the backward Monte-Carlo
difference equation of motion, Eq.~(\ref{eq:backward_motion}), leads to
improved time convergence of the backward solution,
Eq.~(\ref{eq:backward_solution}). Finally, a summary of the main findings of
this article follows in section~\ref{sec:summary}.

\section{Sources of error}

\label{sec:errors}

The greatest disadvantage of the Monte-Carlo method is the unavoidable
statistical error caused by the use of random numbers. With the conventional,
forward Monte-Carlo method, the $\delta$-functions in the solution [see
Eq.~(\ref{eq:forward_solution})] make binning of the particles necessary.
With the number of particles in each bin proportional to $N/N_{\mathrm{bin}}$,
the statistical error becomes $\mathcal{O}((N/N_{\mathrm{bin}})^{-1/2})$ with
obvious notations. Backing out the solution by linear interpolation between
the bins introduces an error $\mathcal{O}(N_{\mathrm{bin}}^{-2})$. \\

As was mentioned previously, the backward method was developed for cases where
the solution is of interest only in a point or small part of phase space.
In such situations it reduces the statistical error to $\mathcal{O}(N^{-1/2})$.
Even in cases where the backward method is used to calculate global solutions,
and the statistical error becomes $\mathcal{O}((N/N_{\mathrm{bin}})^{-1/2})$,
it does not have a finite-bin-size error because the absence of $\delta$-%
functions in the solution [see Eq.~(\ref{eq:backward_solution})] makes binning
superfluous. \\

The finite time step, $\Delta{}t$, also introduces an error term. An analysis
of this time-step error requires some care; some mathematical subtleties are
involved. There are two ways through which the time-step error can enter the
numerical solution. \\

First, directly through the numerical approximation of the Feynman-Kac formula
itself. I.e, even if we knew the \emph{exact} stochastic variables $X_{i}(0)$,
the backward solution, Eq.~(\ref{eq:backward_solution}), would have an error
due to the finite time step. To estimate the magnitude of this error, we need
to go back to the derivation of It\^o's formula in Appendix~A of
Ref.~\cite{bmc} where terms of order $\mathcal{O}(\Delta{}t^{3/2})$ and higher
were neglected. However, in the same way that the
$\mathcal{O}(\Delta{}t^{1/2})$-term, which was kept, becomes zero after
averaging, so does the $\mathcal{O}(\Delta{}t^{3/2})$-term. So, even if the
$X_{i}^{\<}(0)$ were known to all orders, $X_{i}^{\<}(0) \equiv X_{i}(0)$,
the backward solution Eq.~(\ref{eq:backward_solution}), would have an error
$\mathcal{O}(\Delta{}t^{2})$. The error of the forward solution,
Eq.~(\ref{eq:backward_solution}) would of course be of the same order. \\

The other way in which the finite time step introduces an error into the
solution is through an error in $X_{i}^{\<}(0)$. Again, going back to
Ref.~\cite{bmc}, we see that the drift term, $\mu \Delta{}t$, on the right-hand
side of the backward Monte-Carlo equation of motion,
Eq.~(\ref{eq:backward_motion}), has an error $\mathcal{O}(\Delta{}t^{3/2})$,
and the diffusive term, $\zeta \sigma \sqrt{\Delta{}t}$, has an error
$\mathcal{O}(\Delta{}t)$. The different character of these two error terms
should be noted; the $\mathcal{O}(\Delta{}t^{3/2})$-term is deterministic
whereas the $\mathcal{O}(\Delta{}t)$-term is stochastic.
It can be shown (see e.g.~Ref.~\cite{milshtein2}) that after averaging over all
the $X_{i}^{\<}(0)$, the approximation Eq.~(\ref{eq:backward_motion})
introduces an $\mathcal{O}(\Delta{}t)$ error into the backward solution.
Again, the situation is equivalent for the forward method. \\

In section~\ref{sec:higher} we will discuss how the time-step error
of the backward solution can be reduced. But before that, we will study the
backward and forward solutions to a simple model equation in the next section.

\section{Validation and benchmarking}

\label{sec:model}

To validate and benchmark the backward method and to compare it to the
conventional forward method, we will solve a simple model equation that has
an analytic solution. We have chosen the Lorentz equation that has been used
to model pitch-angle scattering in plasmas~\cite{lorentz}:
\begin{equation}
\frac{\partial f}{\partial t} =
\frac{\partial}{\partial x} \, D \,
\frac{\partial f}{\partial x} \ ,
\quad -1 \leq x \leq 1 \ , \quad 0 \leq t \leq T \ , \label{eq:lorentz}
\end{equation}
where the diffusion coefficient $D = (1 + x) (1 - x)$ and the initial
condition is $f(x,0) = \sum_{\ell = 0}^{\infty} (\ell + 1/2) P_{\ell}(x_{0})
P_{\ell}(x) e^{- \ell (\ell + 1) T_{0}}$, where $P_{\ell}(x)$ is the
Legendre polynomial. The analytic solution is
$f(x,T) = \sum_{\ell = 0}^{\infty} (\ell + 1/2) P_{\ell}(x_{0})
P_{\ell}(x) e^{- \ell (\ell + 1) (T_{0} + T)}$. \\

In the following we will use the parameters: $x_{0} = -0.9$, $T_{0} = 0.1$,
and $T = 0.1$; while $N$, $N_{\mathrm{bin}}$, and $\Delta{}t$ will be varied
so that their impact on the error can be studied. Before we get into a
detailed analysis of accuracy and convergence, we show the solutions (with
$N = 2 \times 10^{5}$, $N_{\mathrm{bin}} = 20 $, and $\Delta{}t= 10^{-2}$) to
the Lorentz equation~(\ref{eq:lorentz}) in Fig.~\ref{fig:bulk+tail} below.
\begin{figure}[!t]
\begin{center}
\parbox[c]{.49\textwidth}{%
\includegraphics[scale=0.55,angle=90]{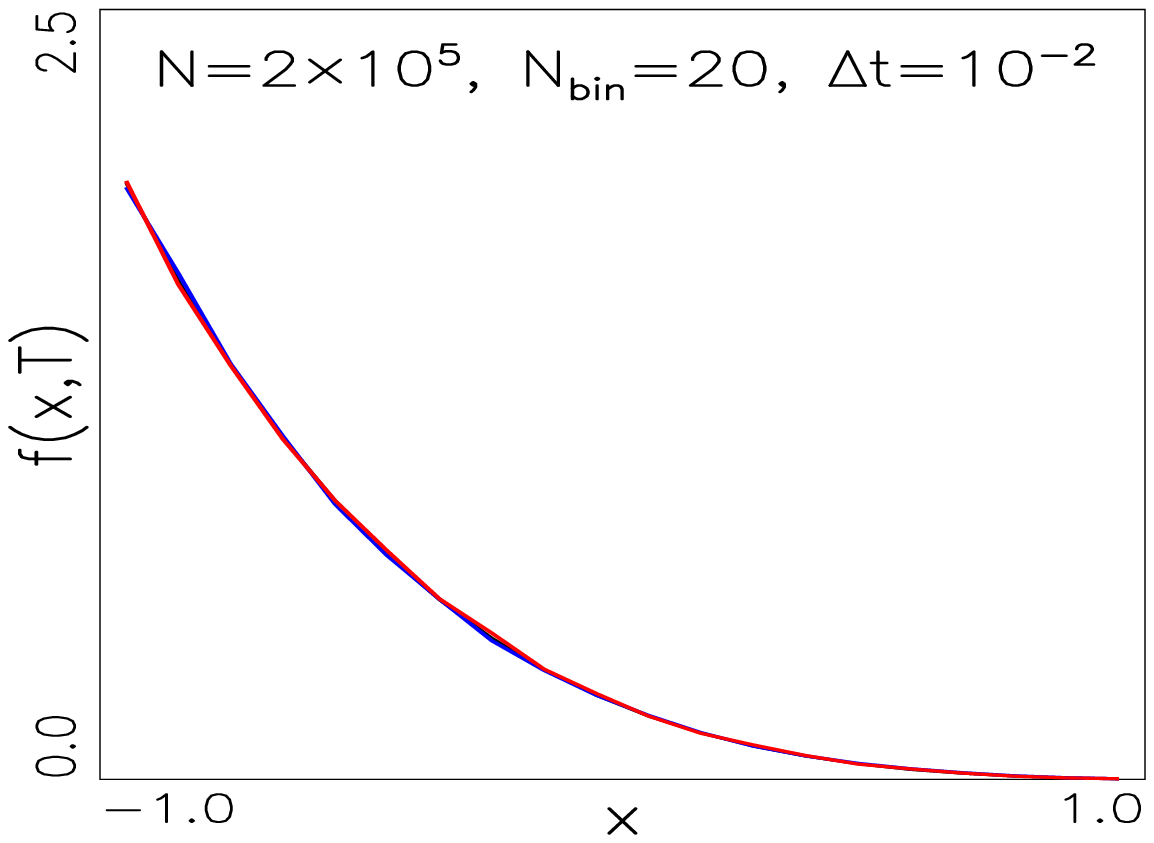}}
\parbox[c]{.49\textwidth}{%
\includegraphics[scale=0.55,angle=90]{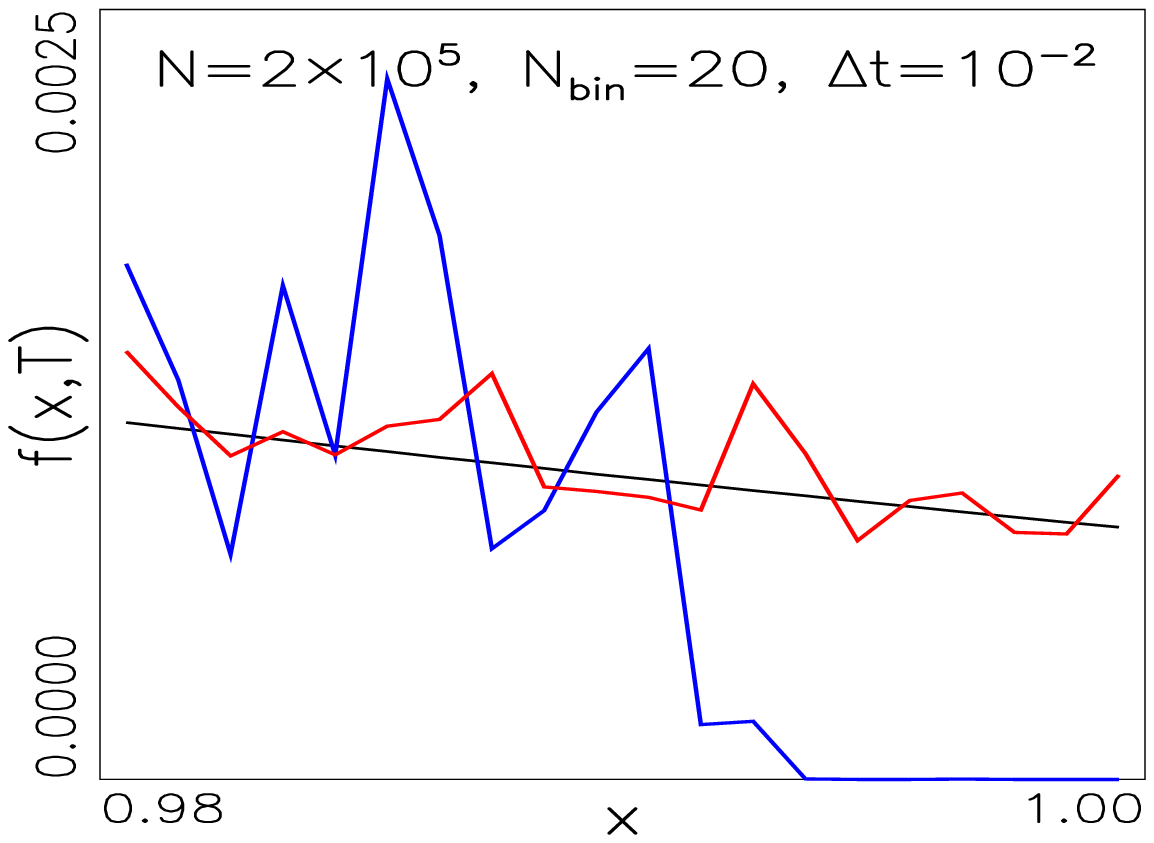}}
\end{center}
\caption{Black line is analytic solution, red line is backward, and blue line
is forward Monte-Carlo solution.
Global solutions are to the left and tail solutions to the right.}
\label{fig:bulk+tail}
\end{figure}
Solving for $f(x,T)$ over the whole interval $-1 \leq x \leq 1$, both the
backward, $f_{\!\<}(x,T)$, and the forward, $f_{\!\>}(x,T)$ Monte-Carlo
solutions are indistinguishable from the exact solution. However, solving only
in the subinterval $0.98 \leq x \leq 1$, the backward solution is a much better
approximation of the exact solution due to its much smaller statistical error.
Note especially that the forward solution drops to zero for $x \gtrsim 0.993$
because not a single particle finds its way into the last few bins
[compare also with the sketch in Fig.~\ref{fig:forward}]. The ability to focus
in on a small subinterval and calculate efficiently (i.e.~without wasting the
vast majority of the particles) the local solution is one of the main
strengths of the backward Monte-Carlo method. \\

For a detailed error analysis we need a well-defined measure of the error.
We will use the local relative error
$\eps{}_{\!\=} = | f_{\!\=}(x_{0},T) - f(x_{0},T) | \, / \, f(x_{0},T)$.
In the following we will plot the logarithm of the error against the logarithm
of one of the parameters $N$, $N_{\mathrm{bin}}$, and $\Delta{}t$, while
keeping the other two parameters fixed. First we compare the statistical error
of the forward and the backward solutions. To single out the statistical error
from the finite bin-size and time-step errors we make the latter two small by
choosing $N_{\mathrm{bin}} = 20$, and $\Delta{}t = 10^{-4}$.
With $N$ spanning six decades ($N = 1,2,5,10,20,50,\ldots,10^{6}$), we plot
the error in the left frame of Fig.~\ref{fig:N-scaling} below.
\begin{figure}[!t]
\begin{center}
\parbox[c]{.49\textwidth}{%
\includegraphics[scale=0.55,angle=90]{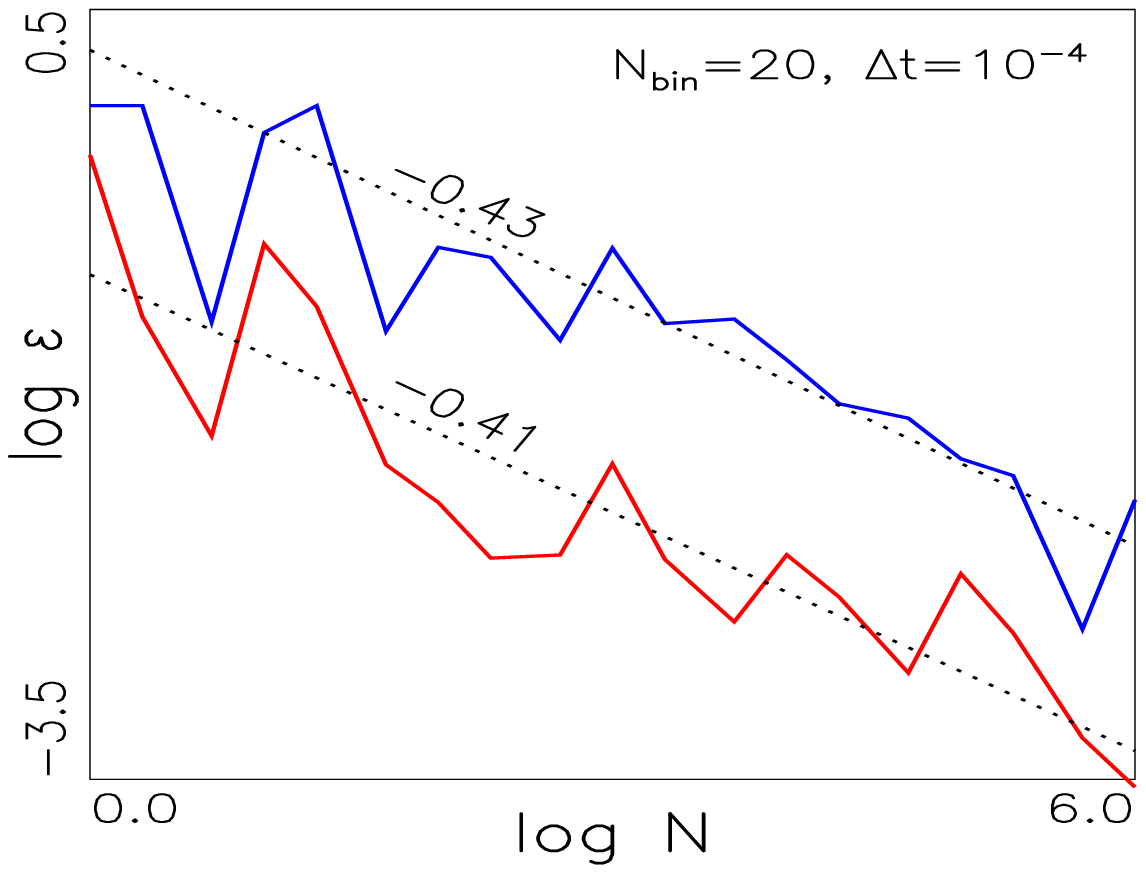}}
\parbox[c]{.49\textwidth}{%
\includegraphics[scale=0.55,angle=90]{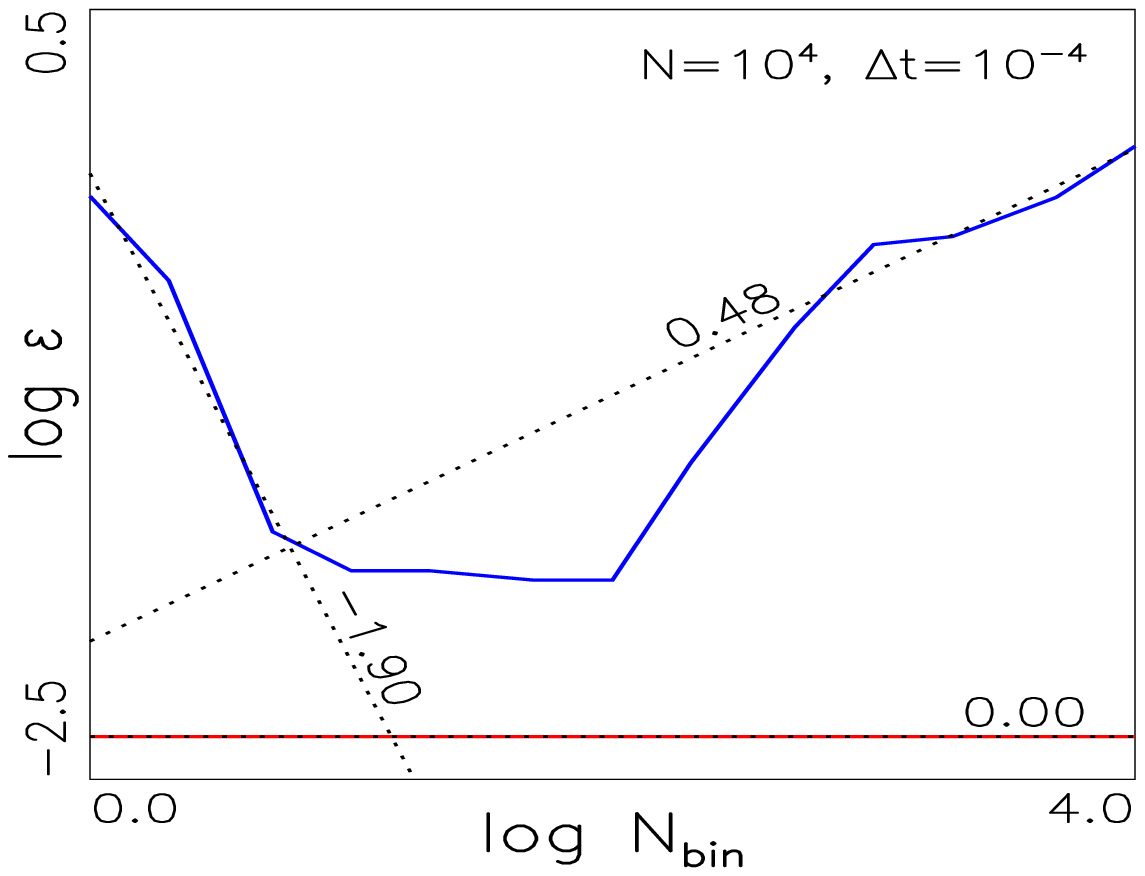}}
\end{center}
\caption{Scaling of the error in the forward (blue) and backward (red)
solutions. The error as a function of the number of particles ($N$) is
plotted in the left frame, and as a function of the number of bins
($N_{\mathrm{bin}}$) to the right.}
\label{fig:N-scaling}
\end{figure}
The dotted lines are least-square fits whose slopes approximate the exponents
that determine the scaling of the error. With the forward method, the bins
adjacent to $x_{0}$ are empty for small $N$ ($N = 1,2,20$). As a result,
$\eps{}_{\!\>} = 1$ for these values of $N$, and the slope becomes flat. The
low-$N$ data points have thus been excluded from the least-square fit to
$\log \eps{}_{\!\>}$. As can be seen, the statistical error
[$\mathcal{O}(N^{-1/2})$ with $N_{\mathrm{bin}}$ fixed] then scales as
predicted. The exponents -0.43 and -0.41 are as close to the theoretical value
of -1/2 as could be expected, given the fact that the other small but non-zero
error terms always tend to flatten the slope. The forward and backward methods
thus converge at the same rate as $N$ is increased, but the forward method
needs $N_{\mathrm{bin}}$ times more particles to achieve the same accuracy
as the backward method. This can easily translate into orders of magnitude in
terms of execution time. \\

Next, we study the bin-size error (and the statistical error) by varying
$N_{\mathrm{bin}}$ over four decades ($N_{\mathrm{bin}} = 1,2,5,\ldots,10^{4}$)
while keeping fixed $N = 10^{4}$ and $\Delta{}t = 10^{-4}$. The result is shown
in the right frame of Fig.~\ref{fig:N-scaling} above. The first observation is
that the backward solution is completely unaffected by the bin size;
$\eps{}_{\!\<}$ is constant and the slope is zero. The error in the forward
solution, $\eps{}_{\!\>}$, exhibits a more interesting behavior with different
scalings for small and large $N_{\mathrm{bin}}$. With few bins there are many
particles in each bin, so the statistical error is small. The bin-size error,
$\mathcal{O}(N_{\mathrm{bin}}^{-2})$, caused by backing out the solution by
linear interpolation between these huge bins, however, becomes large. For small
$N_{\mathrm{bin}}$, the slope of $\log \eps{}_{\!\>}$ is -1.90, in excellent
agreement with the theoretical value of -2. With many bins, the bin-size
error becomes small, but with few particles in each bin, the statistical error
[$\mathcal{O}(N_{\mathrm{bin}}^{1/2})$ with $N$ fixed] becomes dominant as is
evidenced by the slope ($0.48 \lesssim 1/2$) for large $N_{\mathrm{bin}}$.
The shape of the $\log \eps{}_{\!\>}$ graph illustrates the conflict between
resolution and statistics that is inherent for the forward Monte-Carlo method.
The backward solution, however, can be calculated in two points arbitrarily
close without affecting the statistical error. This is yet another of the main
strengths of the backward method. \\

Turning now to the time-step error, we fix $N = 2 \times 10^{5}$ and
$N_{\mathrm{bin}} = 20$ and vary
$\Delta{}t = 10^{-5}, 2 \times 10^{-5}, 5 \times 10^{-5}, \ldots, 10^{-1}$.
As can be seen in Fig.~\ref{fig:t-low} below,
\begin{figure}[!t]
\begin{center}
\includegraphics[scale=0.55,angle=90]{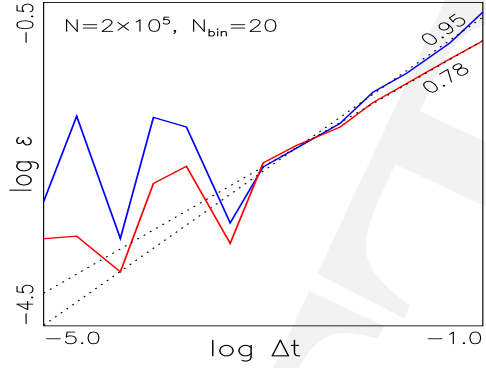}
\end{center}
\caption{The error as a function of the time step ($\Delta{}t$) with the
forward (blue) and backward (red) methods.}
\label{fig:t-low}
\end{figure}
both the forward and the backward methods converge as predicted
(the respective slopes are $\lesssim 1$) down to a time step
$\Delta{}t \approx 10^{-3}$. At this point the statistical error becomes
dominant. \\[5mm]

\section{A higher-order backward method}

\label{sec:higher}

\begin{figure}[!b]
\begin{center}
\parbox[c]{.49\textwidth}{%
\includegraphics[scale=0.55,angle=90]{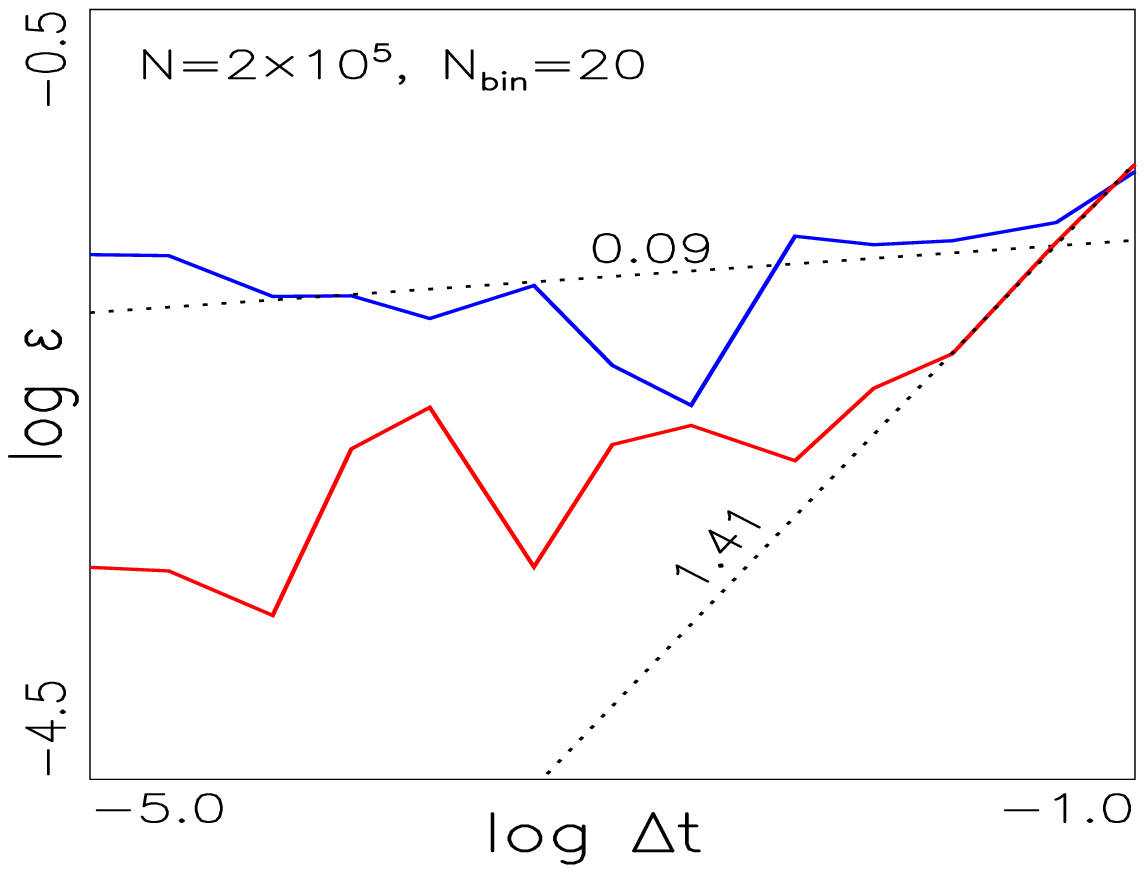}}
\parbox[c]{.49\textwidth}{%
\includegraphics[scale=0.55,angle=90]{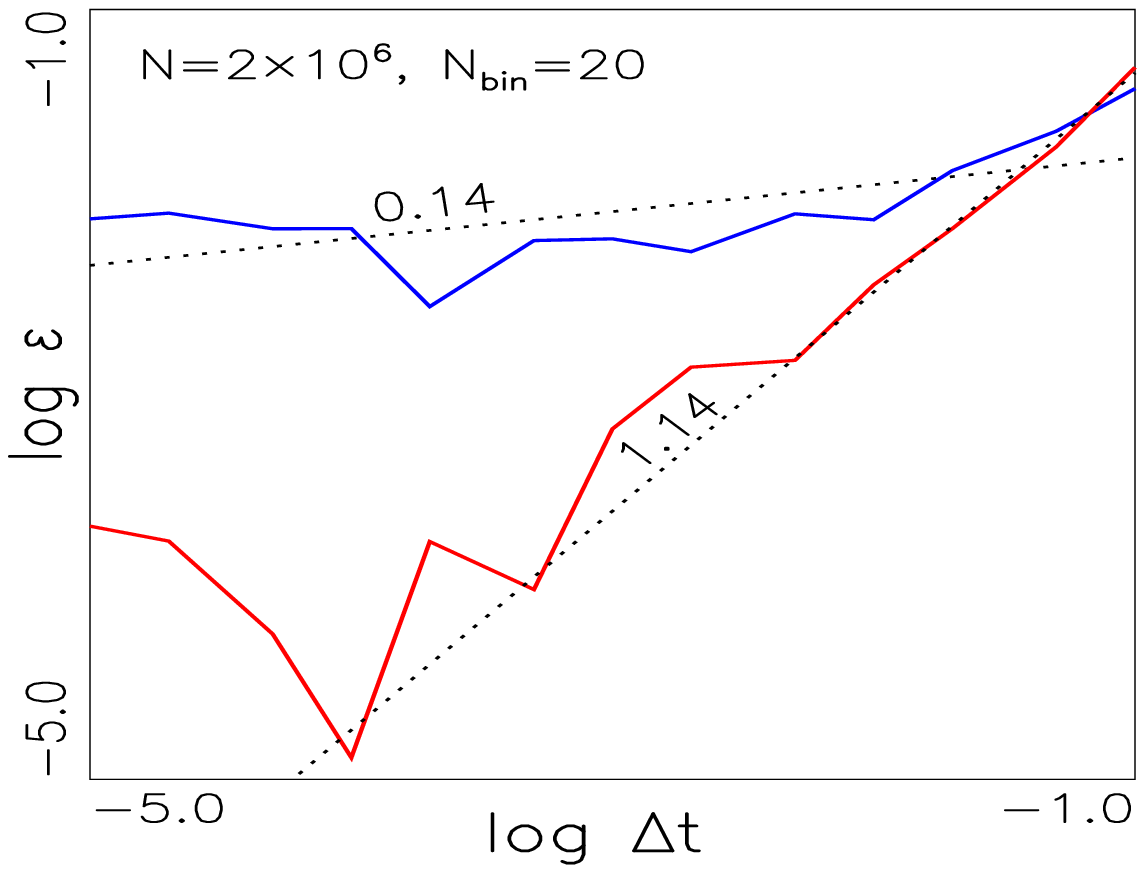}}
\end{center}
\caption{The error as a function of the time step ($\Delta{}t$) with the
higher-order forward (blue) and backward (red) methods with two hundred
thousand particles (left frame) and two million particles (right frame).}
\label{fig:t-high}
\end{figure}

As was discussed in section~\ref{sec:errors}, the dominant
$\mathcal{O}(\Delta{}t)$ time-step error in the Monte-Carlo equations of
motion, Eqs.~(\ref{eq:forward_motion}) and~(\ref{eq:backward_motion}),
comes from the diffusive term $\zeta \sigma \sqrt{\Delta{}t}$. There is thus
reason to suspect that a higher-order diffusive term might result in better
overall time convergence, by reducing the time-step error to
$\mathcal{O}(\Delta{}t^{3/2})$. In appendix~\ref{app:trapezoidal} it is
shown that the next-higher-order Monte-Carlo equation of motion is
\begin{equation}
X_{i}^{\=}(t \pm \Delta{}t) = X_{i}^{\=}(t) +
\half (1 + \zeta^{2}) \, \mu \Delta{}t +
\zeta \sigma \sqrt{\Delta{}t} + \mathcal{O}({\Delta{}t}^{3/2})
\ . \label{eq:backward_motion_higher}
\end{equation}
In the left frame of Fig.~\ref{fig:t-high} below the time-convergence study
presented in Fig.~\ref{fig:t-low} above is repeated using the higher-order
approximation of Eq.~(\ref{eq:backward_motion_higher}). In the right frame
exactly the same experiment is repeated with ten times more particles
to further reduce the statistical error.
Somewhat surprisingly, the forward solution does not converge at all.
The backward solution, however, exhibits the faster time convergence we had
hoped for (the slope is $\lesssim 3/2$). \\

Milshtein, who was not explicitly concerned with solving parabolic equations,
did investigate the scaling of the error of expectation values similar to the
ones of interest here~\cite{milshtein2}. His results are consistent with the
$\mathcal{O}(\Delta{}t)$ error of the forward and backward lower-order
solutions. However, for the higher-order methods introduced in this section,
Milshtein's results would still indicate a $\mathcal{O}(\Delta{}t)$ error for
both the forward and backward solutions, in clear disagreement with the
scalings of Fig.~\ref{fig:t-high}. We hope to be able to resolve the
discrepancies in future work.

\section{Summary}

We have shown that the backward Monte-Carlo method works as expected, i.e.~%
it dramatically reduces the statistical error in situations where the solution
is sought only in a small part of phase space. Furthermore, even in cases where
we solve for the global solution, the backward method removes the conflict
between resolution and statistics. This has great practical significance e.g.~%
when the gradient of the solution is of interest. \\

We have also shown how a remarkably simple modification of the backward
Monte-Carlo equation of motion leads to improved time convergence.

\label{sec:summary}

\begin{ack}

This research was sponsored by the Oak Ridge National Laboratory,\linebreak
managed by UT--Battelle, LLC, for the U.S.~Department of Energy under\linebreak
contract DE--AC05--00OR22725. Research was supported in part by an\linebreak
appointment to  the ORNL Postdoctoral Research Associates Program,\linebreak
administered jointly by Oak Ridge National Laboratory and the Oak Ridge
Institute for Science and Education. \\

The author wishes to thank his colleagues in the Fusion Energy Division
\linebreak
Radio\-frequency Theory Group (Don Batchelor, Lee Berry, Mark Carter, and
Fred Jaeger) for helpful comments during the work on this article.

\end{ack}

\newpage

\appendix

\section{A new perspective on the forward Monte-Carlo method}

\label{app:fmc}

The aim of this appendix is to derive the conventional, forward, weighted
Monte-Carlo method in the same manner as the backward method was
derived~\cite{bmc}. The starting point is the \emph{forward} SDE:
\begin{equation}
dX^{\>} = \mu \, dt + \sigma \, dW \ . \label{eq:forward_sde}
\end{equation}
Solving this SDE is trivial; following exactly the same procedure as for
the backward SDE, we get
\begin{equation}
X^{\>}(t + \Delta{}t) = X^{\>}(t) + \mu \Delta{}t +
\zeta \sigma \sqrt{\Delta{}t} + \mathcal{O}(\Delta{}t) \ , \label{eq:proto0}
\end{equation}
where $\zeta$ is a zero-mean, unit-variance Gaussian random number,
$\zeta \in N(0,1)$. When the SDE is backward, the obvious initial condition is
$X^{\<}(t\!=\!T) = x$. For the forward SDE~(\ref{eq:forward_sde}), we
simply do not know which initial condition to impose; we will have to leave
$X^{\>}(t\!=\!0)$ undefined for now. \\

As before~\cite{bmc}, we use the It\^o formula to differentiate
$f(X(t),t)$ and obtain
\begin{equation}
f(X^{\>}(T),T) = f(X^{\>}(0),0) + \int_{0}^{T} 2 f_{t} \, dt +
\int_{0}^{T} \sigma f_{x} \, dW \ ,
\label{eq:proto1}
\end{equation}
where again we have identified $\mu = \partial D / \partial x$,
$\sigma = \sqrt{2 D}$, and used\linebreak
$f_{t} + \mu f_{x} + \half \sigma^{2} f_{xx} = $
[Eq.~(\ref{eq:canonical_diffusion})] $ = 2  f_{t}$. Using the macroscopic
initial condition, $f(x,0) = \Phi(x)$, we get
\begin{equation}
f(X^{\>}(T),T) = \Phi{}(X^{\>}(0)) - \int_{0}^{T} \sigma f_{x} \, dW \ .
\label{eq:proto2}
\end{equation}
At this stage of the derivation of the \emph{backward} method, we used the
microscopic initial condition and found that the equivalent of the LHS of
Eq.~(\ref{eq:proto2}) was in fact a solution to
Eq.~(\ref{eq:canonical_diffusion}): $f(X^{\<}(t\!=\!T),T) =$
[$X^{\<}(t\!=\!T) = x$] $ = f(x,T)$. But here, in the forward derivation,
we cannot do that because $X^{\>}(T) \neq x$! Going forward in time,
$X^{\>}(T)$ is an unknown; we could use the Monte-Carlo equation
of motion, Eq.~(\ref{eq:proto0}), to find $X^{\>}(T)$, but we do not have an
initial condition! So, is there any way to get $f(x,T)$ from
Eq.(\ref{eq:forward_sde})? The answer is yes; there is a (rather contrived)
way. \\

If we let $f(x,T)$ be a distribution, we can write it as
\begin{equation}
f(x,T) = E[ \, f(X^{\>}(T),T) \times \delta{}(x - X^{\>}(T)) \, ] \ .
\label{eq:proto3}
\end{equation}
Substituting Eq.~(\ref{eq:proto2}) into Eq.~(\ref{eq:proto3}), we get:
\begin{equation}
f(x,T) = E[ \, \Phi{}(X^{\>}(0)) \times \delta{}(x - X^{\>}(T)) \, ] \ .
\label{eq:proto4}
\end{equation}
The numerical approximation of this expectation value is the forward weighted
Monte-Carlo solution:
\begin{equation}
f_{\!\>}(x,T) = N^{-1} \sum_{i=1}^{N} \Phi{}(X_{i}^{\>}(0)) \times
\delta (x - X_{i}^{\>}(T)) \ . \label{eq:proto5}
\end{equation}
Now, we still need to know what the microscopic initial condition on
$X_{i}^{\>}(0)$ should be. The best we can do is to make sure that
$f_{\!\>}(x,T\!=\!0) \approx \Phi{}(x)$. Note that ``approximately equal'' must
be given a very liberal definition because $f_{\!\>}(x,T\!=\!0)$ is a jagged
sum of $\delta$-functions, whereas $\Phi{}(x)$ is in general a smooth function.
The stochastic variables $X_{i}^{\>}(0)$ can now be found by following the
stochastic trajectories given by the forward Monte-Carlo equation of motion,
Eq.~(\ref{eq:proto0}), with the ``fuzzy initial condition'':
\begin{equation}
X_{i}^{\>}(0) \in
\{X_{i}^{\>}(0) \ | \ f_{\!\>}(x,T\!=\!0) \approx \Phi{}(x)\} \ .
\label{eq:fuzzy}
\end{equation}

\section{A higher-order Monte-Carlo equation of motion}

\label{app:trapezoidal}

The diffusive term $\zeta \sigma{}(X(t)) \sqrt{\Delta{}t}$ of
the Monte-Carlo equations of motion,
Eqs.~(\ref{eq:forward_motion}) and~(\ref{eq:backward_motion}),
is an approximation of the It\^o integral
\begin{equation}
\int_{t}^{t + \Delta{}t}\!\sigma{}(X(s)) \, dW(s) =
\sigma{}(X(t)) \int_{t}^{t + \Delta{}t}\!dW(s) +
\mathcal{O}(\Delta{}t) \ ,
\label{eq:sigma_integral_low}
\end{equation}
where $X(s)$ is the stochastic process that solves the SDE
\begin{equation}
dX = \mu{}(X(s)) \, ds + \sigma{}(X(s)) \, dW(s) \ , \label{eq:sde}
\end{equation}
where $W$ is a Wiener process, and we impose the initial condition
$X(t) = x_{t}$.
The approximation Eq.~(\ref{eq:sigma_integral_low}) replaces $\sigma{}(x)$
with its zero-order Taylor expansion, $\sigma_{0}{}(x) = \sigma{}(x_{t})$.
To reduce the error in Eq.~(\ref{eq:backward_motion}) to
$\mathcal{O}({\Delta{}t}^{3/2})$ we Taylor expand $\sigma{}(x)$ to first order,
\begin{equation}
\sigma_{1}{}(x) = \sigma{}(x_{t}) + \sigma{}'(x_{t}) (x - x_{t}) \ ,
\label{eq:sigma_one}
\end{equation}
and approximate the solution to Eq.~(\ref{eq:sde}) with
\begin{equation}
X(s) = x_{t} + \mu{}(x_{t}) (s - t) + \sigma{}(x_{t}) [W(s) - W(t)] \ .
\label{eq:X_approx}
\end{equation}
The It\^o integral over $\sigma$ becomes:
\begin{eqnarray}
\int_{t}^{t + \Delta{}t}\!&\sigma&{}(X(s)) \, dW(s) =
\sigma(x_{t}) \int_{t}^{t + \Delta{}t}\!dW(s) + \nonumber \\
&\sigma&{}'(x_{t}) \int_{t}^{t + \Delta{}t}\!\left \{ \mu(x_{t}) (s - t) +
\sigma(x_{t}) [W(s) - W(t)] \right \} dW(s) +
\mathcal{O}(\Delta{}t^{3/2}) = \nonumber \\
&\sigma&(x_{t}) \int_{t}^{t + \Delta{}t}\!dW(s) +
\sigma{}(x_{t}) \, \sigma{}'(x_{t})
\int_{t}^{t + \Delta{}t} [W(s) - W(t)] \, dW(s) +
\mathcal{O}(\Delta{}t^{3/2}) \ , \nonumber \\
{} \label{eq:sigma_integral_high}
\end{eqnarray}
where $\mu(x_{t}) \int_{t}^{t + \Delta{}t} (s - t) \, dW(s)$ is
$\mathcal{O}(\Delta{}t^{3/2})$ and hence was neglected. The first term on the
RHS of Eq.~(\ref{eq:sigma_integral_high}) is just the low-order approximation
Eq.~(\ref{eq:sigma_integral_low}). The second term is the next-order
correction. To calculate this correction term, we first need some
intermediate results. \\

Following Bj\"ork~\cite{bjork}, we introduce the stochastic variables
\begin{equation}
I_{n} = \sum_{j=0}^{n-1} W(t_{k}) [W(t_{k+1}) - W(t_{k})]  \ , \label{eq:I_n}
\end{equation}
and
\begin{equation}
B_{n} = \sum_{j=0}^{n-1} W(t_{k+1}) [W(t_{k+1}) - W(t_{k})]  \ , \label{eq:B_n}
\end{equation}
where $t_{j} = j \, t / n  \, , \ j = 0, 1, \ldots, n-1$. It trivially
follows that
\begin{equation}
B_{n} + I_{n} = \sum_{j=0}^{n-1} [W(t_{k+1}) + W(t_{k})]
[W(t_{k+1}) - W(t_{k})] = W^{2}(t) \ , \label{eq:B+I}
\end{equation}
and
\begin{equation}
B_{n} - I_{n} = \sum_{j=0}^{n-1} [W(t_{k+1}) - W(t_{k})]^{2} = S_{n}(t) \ .
\label{eq:B-I}
\end{equation}
In Appendix~A of Ref.~\cite{bmc} we showed that
\begin{equation}
\lim_{n \rightarrow \infty} S_{n}(t) = t \ .
\label{eq:limit}
\end{equation}
Subtracting Eq.~(\ref{eq:B-I}) from Eq.~(\ref{eq:B+I}) and taking the limit
$n \rightarrow \infty$, we obtain
\begin{equation}
\lim_{n \rightarrow \infty} I_{n}(t) = \int_{0}^{t} W(s) \, dW(s) =
\half [W^{2}(t) - t] \ . \label{eq:W_integral}
\end{equation}

Applying Eq.~(\ref{eq:W_integral}) to Eq.~(\ref{eq:sigma_integral_high})
we get
\begin{eqnarray}
\int_{t}^{t + \Delta{}t}\!&\sigma&{}(X(s)) \, dW(s) =
\sigma{}(X(t)) \, [W(t + \Delta{}t) - W(t)] + \nonumber \\
&\half& \sigma{}(X(t)) \, \sigma{}'(X(t)) \,
\{[W(t + \Delta{}t) - W(t)]^{2} - \Delta{}t \} +
\mathcal{O}(\Delta{}t^{3/2}) \ . \nonumber \\
{} \label{eq:sigma_integral_high_eval}
\end{eqnarray}
In the absence of advection, as is the case for
Eq.~(\ref{eq:canonical_diffusion}), $\sigma \sigma{}' = \mu$, and the
higher-order Monte-Carlo equations of motion become
\begin{eqnarray}
X_{i}^{\=}(t \pm \Delta{}t) &=&
X_{i}^{\=}(t) + \mu \Delta{}t + \zeta \sigma \sqrt{\Delta{}t} +
\half (\zeta^{2} - 1) \mu \Delta{}t + \mathcal{O}({\Delta{}t}^{3/2}) =
\nonumber \\
&{}& X_{i}^{\=}(t) + \half (1 + \zeta^{2}) \mu \Delta{}t +
\zeta \sigma \sqrt{\Delta{}t} + \mathcal{O}({\Delta{}t}^{3/2}) \ ,
\label{eq:backward_motion_higher_long}
\end{eqnarray}
where we have used the fact that
$W(t + \Delta{}t) - W(t) \in N(0,\sqrt{\Delta t})$;
$\zeta$ is a zero-mean, unit-variance Gaussian random number;
$\zeta \in N(0,1)$; and $\mu = \mu{}(X(t))$ and $\sigma = \sigma{}(X(t))$. \\

In Ref.~\cite{milshtein1} Milshtein derived
Eq.~(\ref{eq:sigma_integral_high_eval}) by introducing
\begin{equation}
X^{*}(t + \Delta{}t) = x_{t} + c_{1} \Delta{}t +
c_{2} [W(t + \Delta{}t) - W(t)] + c_{3} [W(t + \Delta{}t) - W(t)]^{2}  \ ,
\nonumber
\end{equation}
and finding the $c_{1}$, $c_{2}$ and $c_{3}$ that minimize
$E[(X^{*}(t + \Delta{}t) - X(t + \Delta{}t))^{2}]$.

\end{document}